\documentclass[preprint,12pt]{elsarticle}
\usepackage{amssymb}
\usepackage[nodots]{numcompress}
\journal{Mathematics and Computers in Simulation}
\usepackage[
   paper=a4paper,
   tmargin=4cm,           
   bmargin=3cm,           
   lmargin=2.5cm,             
   rmargin=2.5cm]{geometry}   
\usepackage[english]{babel}
\usepackage[latin1]{inputenc}
\usepackage{amssymb}
\usepackage{amsmath}
\usepackage{pb-diagram}
\usepackage{graphicx}
\usepackage{color}
\usepackage{fancyhdr}

\usepackage[T1]{fontenc}        
\usepackage{url}                
\usepackage{vmargin}            

\usepackage{lineno}
\usepackage{setspace}                      
\newcommand{\etal}{\it{et \, al.}}

\newcommand{\bo}[1]{\mathbf{#1}}
\def\indic{\hbox{1\kern-.24em\hbox{I}}}      

 \newcommand{\esp}{\mathbb{E}}
 
\newcommand{\w}{w} 

\newcommand{\x}{x}
\newcommand{\X}{X}

\renewcommand{\baselinestretch}{1.2}
\newtheorem{prop}{Proposition}[section]
\newtheorem{defi}{Definition}[section]

\newtheorem{preuve}{Proof}[section]
\newtheorem{theo}{Theorem}[section]
\newtheorem{rem}{Remark}[section]

\catcode`\"=\active
\catcode`\"=\active
\def"{\og\ignorespaces}
\def"{{\fg}}

\begin{document}
\begin{frontmatter}
\title{Derivative-based global sensitivity measures: general links with Sobol' indices and numerical tests}

\author{M. Lamboni$^{ab}$, B. Iooss$^c$\footnote{Corresponding author: bertrand.iooss@edf.fr, Phone: +33 1877969, Fax: +33 130878213}, A.-L. Popelin$^c$, F. Gamboa$^b$}

\address{$^a$ Universit\'e Paris Descartes, 45 rue des saints P\`eres, F-75006, France \\
 	       $^b$ IMT, F-31062, France \\
         $^c$ EDF R\&D, 6 quai Watier, F-78401, France} 


\begin{abstract}
The estimation of variance-based importance measures (called Sobol' indices) of the input variables of a numerical model can require a large number of model evaluations.
It turns to be unacceptable for high-dimensional model involving a large number of input variables (typically more than ten).
Recently, Sobol and Kucherenko have proposed the Derivative-based Global Sensitivity Measures (DGSM), defined as the integral of the squared derivatives of the model output, showing that it can help to solve the problem of dimensionality in some cases.
We provide a general inequality link between DGSM and total Sobol' indices for input variables belonging to the class of Boltzmann probability measures, thus extending the previous results of Sobol and Kucherenko for uniform and normal measures.
The special case of log-concave measures is also described.
This link provides a DGSM-based maximal bound for the total Sobol indices.
Numerical tests show the performance of the bound and its usefulness in practice.
\end{abstract}


\end{frontmatter}

\setpagewiselinenumbers
\modulolinenumbers[1]
\linenumbers
\doublespacing

\section{Introduction}

With the advent of computing technology and numerical methods, computer models are now widely used to make predictions on little-known physical phenomena, to solve optimization problems or to perform sensitivity studies.
These complex models often include hundreds or thousands uncertain inputs, whose uncertainties can strongly impact the model outputs (De Rocquigny $\etal$ \cite{derdev08}, Kleijnen \cite{kle08}, Patelli $\etal$ \cite{patelli10}). In fact, it is well known that, in many cases, only a small number of input variables really act in the model (Saltelli $\etal$ \cite{salcha00}). This number is referred to the notion of the effective dimension of a function (Caflish $\etal$ \cite{caflish1997}), which is a useful way to deal with the curse of dimensionality in practical applications.

Global Sensitivity Analysis (GSA) methods (Sobol \cite{sobol93}, Saltelli $\etal$ \cite{salcha00}) are used to quantify the influence of model input variables (and their interaction effects) on a model reponse. It is also an objective way to determine the effective dimension by using the model simulations (Kucherenko $\etal$ \cite{Kucherenko2011}).
A first class of GSA methods, called ``screening'' methods, aim at dealing with a large number of input variables (from tens to hundreds).
An example of screening method is the Morris' method (Morris \cite{morris91}), which allows a coarse estimation of the main effects using only a few model evaluations.
While taking into account the interactions between the indices, the basic form of the Morris method did not compute precise sensitivity indices associated to the interactions between inputs.
The second class of GSA methods are the popular quantitative methods, mainly based on the decomposition of the model output variance, which leads to the so-called variance-based methods and Sobol' sensitivity indices.
It allows computing the main and total effects (called first order and total Sobol' indices) of each input variable, as well as interaction effects.
However, for functions with non linear and interaction effects, the estimation procedures become particularly expensive in terms of number of required model evaluations.
Hence, for this kind of model, variance-based methods can only be applied to a limited number of input variables (less than tens).

Recently, Sobol and Kucherenko \cite{sobol09,sobol10} have proposed the so-called Derivative-based Global Sensitivity Measures (DGSM), which can be seen as a kind of generalization of the Morris screening method.
DGSM seem computationally more tractable than variance-based measures, specially for high-dimensional models.
They also theoretically proved an inequality linking DGSM to total Sobol' indices in the case of uniform or Gaussian input variables.

In this paper, we investigate this close relationship between total Sobol' indices and DGSM, by extending this inequality to a large class of Boltzmann probability measures. We also obtain result for the class of log-concave measures. The paper is organized as follows: Section \ref{sec:def} recalls some useful definitions of Sobol' indices and DGSM. Section \ref{sec:link} establishes an inequality between these indices for a large class of Boltzmann (resp. log-concave) probability measures. Section \ref{sec:test} provides some numerical simulations on two test models, illustrating how DGSM can be used in practice. We conclude in Section \ref{sec:con}.

\section{Global sensitivity indices definition} \label{sec:def}

\subsection{Variance-based sensitivity indices}

Let $ Y= f(\bo{\X})$ be a model output with $d$ random input variables $\bo{\X}=(X_1,\ldots,X_d)$. If the input variables are independent (assumption A1) and $\esp\left( f^2(\bo{\X}) \right)< +\infty $ (assumption A2), we have the following unique Hoeffding decomposition (Efron and Stein \cite{efron81}) of $ f(\bo{\X})  $:
\begin{eqnarray}  \label{eq:decsa}
  f(\bo{\X})  & = & f_0 + \sum_j^d f_j(\X_j) + \sum_{i<j}^d f_{ij}(\X_i,\X_j) + \ldots + f_{1\ldots d}(\X_1,\ldots,X_d)   \\
             & = & \sum_{u \subset \{1,2, \ldots d\}} f_u(\X_u),
\end{eqnarray}
where $f_0 = \esp \left[ f(\bo{X})\right]  $ corresponds to the empty subset; $f_j(\X_j) = \esp \left[f(\bo{\X})| \X_j  \right] -f_0 $  and $\displaystyle f_u(\X_u) =  \esp \left[f(\bo{\X})| \X_u  \right] - \sum_{v \subset u } f_v(\X_v) $  for any subset $u \subset \{1,2, \ldots, d \}$ .

By regrouping all the terms in equation (\ref{eq:decsa}) that contain the variable $\X_j$ ($j=1,2, \ldots, d $) in the function called $g(\cdot)$:

\begin{equation}\label{eq:g} 
g(\X_j,\bo{\X}_{\sim j}) = \sum_{u \ni j}f_u(\bo{\X}_u)\;,
\end{equation}

we have the following decomposition:
\begin{equation} \label{eq:decdgsa}
 f(\bo{\X}) = f_0+ g(\X_j,\bo{\X}_{\sim j}) + h(\bo{\X}_{\sim j}),
\end{equation}
where $ \bo{\X}_{\sim j}$ denotes the vector containing all variables except $\X_j$ and $h(\cdot)= f(\cdot)- f_0 - g(\cdot)$. Notice that
this decomposition is also unique under assumptions A1 and A2.
The function $g(\cdot)$, itself, suffices to compute the total sensitivity indices. Indeed, it contains all information relating $f(\bo{X})$ to $\X_j$.

\begin{defi}
 Assume that A1, A2 hold, let $\mu(\bo{X})= \mu(\X_1, \ldots, \X_d)$ be the distribution of the input variables. For any non empty subset $u \subseteq \{1,2, \ldots, d \}$, set first
$$D=\int f^2(\bo{\x}) d\mu(\bo{\x}) -f_0^2 \;,$$
$$ D_u=\int  f^2_u(\bo{\x}_u) d\mu(\bo{\x}_u)\;,$$
\begin{equation} \label{eq:dt0}
D_{u}^{tot} = \int \sum_{v \supseteq u } f_v^2(\bo{\x}_v) d\mu(\bo{\x}_v)\;.
\end{equation}
Further, the first order Sobol sensitivity indices (Sobol \cite{sobol93}) of $\bo{\X}_u$ is
 \begin{equation}\label{eq:si}
S_{u} = \frac{D_u}{D} \;,
\end{equation}
The total sensitivity Sobol index of $\bo{\X}_{u}$ (Homma and Saltelli \cite{homsal96}) is
\begin{equation}\label{eq:sit}
S_{T_u} = \frac{D_{u}^{tot}}{D} \;.
\end{equation}
\end{defi}

The following proposition gives another way to compute the total sensitivity indices.
\begin{prop}
Under assumptions A1 and A2, the total sensitivity indices of
variable $\X_j$ ($j=1,2, \ldots, d $) is obtained by the following formulas:
\begin{eqnarray}\label{eq:dt}
D_{j}^{tot}  & = &  \int  g^2(\x_j,\bo{\x}_{\sim j}) d\mu(\bo{\x})
\end{eqnarray}
and
\begin{eqnarray}\label{eq:dt2}
D_{j}^{tot} & = &  \frac{1}{2}\int \left[f(\bo{\x})-f(\x'_j, \bo{\x}_{\sim j}) \right]^2 d\mu(\bo{\x})
d\mu(\x_j')\;.
\end{eqnarray}
\end{prop}

\begin{preuve}
The first formula is an obvious consequence of equation (\ref{eq:decdgsa}), and it is obtained by using the orthogonality of the summands in equation (\ref{eq:decsa}). Indeed, $\displaystyle D_{j}^{tot} = \int \sum_{v \supseteq j} f_v^2(\bo{\x}_v) d\mu(\bo{\x}_v) = \int \left[\sum_{v\supseteq j} f_v(\bo{\x}_v) \right]^2  d\mu(\bo{\x})=    \int  g^2(\x_j,\bo{\x}_{\sim j}) d\mu(\bo{\x})$.
 The later formula is proved in Sobol \cite{sobol2001}.
\end{preuve}

\subsection{Derivative-based sensitivity indices}

Derivative-based global sensitivity method uses the second moment of model derivatives as importance measure. This method is motivated by the fact that a high value of the derivative of the model output with respect to some input variable means that a big variation of model output is expected for a variation of the variable. This method extends the Morris method (Morris \cite{morris91}). Indeed, it allows to capture any small variation of the model output due to input variables.

DGSM have been first proposed in Sobol and Gresham \cite{sobgre95}. Then, they have been largely studied in Kucherenko $\etal$ \cite{kucherenko09}, Sobol and Kucherenko \cite{sobol09,sobol10} and Patelli $\etal$ \cite{patelli10b}.
From now on, we assume that the function $f$ is differentiable.
Two kind of DGSM are defined below:

\begin{defi}
Assume that A1 holds and that $\displaystyle \frac{\partial f(\bo{\X})}{\partial \x_j}$ is square-integrable (assumption A3). Then, for $ j=1,2, \ldots d $, we define the DGSM indices by:
\begin{eqnarray}\label{eq:dsi}
\nu_{j}  & = &  \esp \left[\left(\frac{\partial f(\bo{\X})}{\partial \x_j}\right)^2\right]   \\ \nonumber
   & = & \int  \left(\frac{\partial f(\bo{\x})}{\partial \x_j}\right)^2 d\mu(\bo{\x}) \;.
\end{eqnarray}
 Let $w(\cdot)$ is be a bounded measurable function. A weighted version of the last indices is:
\begin{equation} \label{eq:dsip}
\tau_{j}  = \int \left(\frac{\partial f(\bo{\x})}{\partial \x_j}\right)^2 w(\x_j) d\mu(\bo{\x}) .
\end{equation}
\end{defi}

\begin{rem} \label{rem:sieqdsi}
Sobol and Kucherenko \cite{sobol10} showed that, for a specific weighting function $\displaystyle w(\x_j) = \frac{1-3\x_j+3 \x_j^2}{6} $ and for a class of linear model with respect to each input variable (following a uniform distribution over $[0,1]$), we have  $ \tau_j = D_{j}^{tot} $.
\end{rem}

 \begin{rem}
By bearing in mind the decomposition in equation (\ref{eq:decdgsa}), we can replace in equations (\ref{eq:dsi} ) and (\ref{eq:dsip}) the function $f(\cdot)$ by $g(\cdot)$.  In general, $g(\cdot)$ is a $d_1$ ($d_1 \leq d$) dimension function, and this can drastically reduce the number of model evaluations for the numerical computation of $\nu$ or $\tau$. Thus, we have:
\begin{eqnarray}\label{eq:dsi2}
\nu_{j}  & = & \int  \left(\frac{\partial g(\bo{\x})}{\partial \x_j}\right)^2 d\mu(\bo{\x}) \;.
\end{eqnarray}
\begin{equation} \label{eq:dsip2}
\tau_{j}  = \int \left(\frac{\partial g(\bo{\x})}{\partial \x_j}\right)^2 w(\x_j) d\mu(\bo{\x}) ,
\end{equation}
\end{rem}
%

\section{Variance-based sensitivity indices vs. derivative-based sensitivity indices} \label{sec:link}

As DGSM estimations need much less model evaluations than total Sobol' indices estimations (Kucherenko $\etal$ \cite{kucherenko09}), it would be interesting to use the DGSM, instead of total Sobol' indices, for factors fixing setting.
A formal link is therefore necessary to provide a mathematical relation between total Sobol' indices and DGSM.
Sobol and Kucherenko \cite{sobol09} have established an inequality linking these two indices for uniform and Gaussian random variables (maximal bound for $S_{T_j}$).
In this section, we extend the inequality for Sobolev' space model whith the marginal distribution of input variables belonging to the class of Boltzmann measure on $\mathbb{R}$ (assumption A4).
A measure $\delta$ on $\mathbb{R}$ is said to be a Boltzmann measure if it is absolutely continuous with respect to the Lebesgue measure and its density $\displaystyle d\delta(\x) = \rho(\x)d\x = c \exp[-v(\x)]d\x$. Here $v(\cdot)$ is a continuous function and $c$ a normalizing constant.
 Many classical continuous probability measures used in practice are Boltzmann measures (see de Rocquigny $\etal$ \cite{derdev08} and Saltelli $\etal$ \cite{salcha00}).

The class of Boltzmann probability measures includes the well known class of log-concave probability measures. In this case, $v(\cdot)$ is a convex function (assumption A5).
In other words, a twice differentiable probability density function $\rho(\x)$ is said to be log-concave if, and only if,
\begin{equation}
\frac{d^2}{dx^2}[\log \rho(\x)] \leq 0 \;.
\end{equation}
Note that the probability measure of uniform density on a finite interval is not continuous on $\mathbb{R}$.
So it cannot be considered in the class of log-concave probability measure, nor in the class of Boltzmann probability measures.

The two following propositions give the formal link between Sobol' indices and derivative-based sensitivity indices.

\begin{theo} \label{prop:link}
Under assumptions A1, A2, A3 and A4, we have:
\begin{equation}\label{eq:ineq}
     D_{j}^{tot} \leq C(\mu_j) \nu_{j} \,
\end{equation}
 with $\displaystyle C(\mu_j) = 4C_1^2$ and $\displaystyle C_1=\sup_{\x \in \mathbb{R}} \frac{\min(F_j(\x), 1-F_j(\x))}{\rho_j(\x)} $ the Cheeger constant, $F_j(\cdot)$ the cumulative probability function of $\X_j$ and $\rho_j(\cdot)$ the density of $\X_j$.
\end{theo}
We recall the four assumptions:
\begin{itemize}
\item A1: independence between inputs $X_1$, $X_2$, \dots, $X_d$,
\item A2: $f \in L^2(\mathbb{R})$,
\item A3: $\displaystyle \frac{\partial f}{\partial x_j} \in L^2(\mathbb{R})$,
\item A4: the distribution of $X_j$ is a Boltzmann  probability measure.
\end{itemize}


 \begin{preuve} \label{proof:link1}
The resulting inequality (\ref{eq:ineq}) is based on a one-dimensional $L^2$-Poincar\'e inequality of the type $\displaystyle \|u\|_{L^2} \leq C  \|\nabla u\|_{L^2}$ for $u$ a Sobolev' space function (see for example \cite{fougeres05}).
It is applied here to the function $g(\cdot)$ (equation (\ref{eq:g}), with $\displaystyle  \int g^2(\x_j,\bo{\x}_{\sim j}) d\mu(\bo{\x}) = D_{j}^{tot}$ (equation (\ref{eq:dt})) and $\displaystyle  \int  \left(\frac{\partial g(\bo{\x})}{\partial \x_j}\right)^2 d\mu(\bo{\x}) = \nu_{j}$ (equation (\ref{eq:dsi2})).
 The constant is obtained in Bobkov \cite{bobkov99}, and Foug\`eres \cite{fougeres05} for the one-dimensional Poincar\'e inequality. A proof of the $d$-dimensional Poincar\'e inequality is given in Bakry $\etal$ \cite{bakry08}.
 \end{preuve}

\begin{theo} \label{prop:link2}
Under assumptions A1, A2, A3 and A5, we have:
\begin{equation}\label{eq:ineq2}
     D_{j}^{tot} \leq  \left[\exp(v(m))\right]^2 \nu_{j} \;,
\end{equation}
 with $\displaystyle C_1 = \frac{\exp(v(m))}{2}$ the Cheeger constant and $m$ the median of the measure $\mu_j$ (such that $\mu(X_j \leq m) = \mu(X_j>m)$).
\end{theo}

We recall the assumption A5: the distribution of $\X_j$ is a log-concave probability measure.

 \begin{preuve}
 See proof \ref{proof:link1}.
 \end{preuve}

Table \ref{tab:logconc} shows Cheeger constant for some log-concave probability distributions that are used in practice for uncertainty and sensitivity analyses. We also give their medians and the functions $v(\cdot)$.
We obtain the same results for the normal distribution ${\mathcal N}(\mu, \sigma^2)$ similar to Sobol and Kucherenko \cite{sobol09} but we prove them in another way (in this case, $v(m) = \log(\sigma)$).
For uniform distribution ${\mathcal U}[a \, b]$, Sobol and Kucherenko \cite{sobol09} obtained via direct integral manipulations the inequality $\displaystyle D_{j}^{tot} \leq  \frac{(b-a)^2}{\pi^2} \nu_{j}$. This relation is the classical Poincar\'e or Writtinger inequality (Ane $\etal$ \cite{ane2000}). 

 \begin{table}[!ht]
\hspace{-1cm}
    \begin{tabular}{p{3.3cm} c c c}
\hline
Distribution & $v(x)$ & $m$ & $C_1$\\
\hline
Normal ${\mathcal N}(\mu, \sigma^2)$ & $\displaystyle \frac{(x-\mu)^2}{2\sigma^2}+\log(\sigma)$ & $\mu$ & $\displaystyle \frac{\sigma}{2}$ \\
\\
Exponential ${\mathcal E}(\lambda)$, $\lambda>0$ & $\lambda x -\log(\lambda)$ & $\displaystyle \frac{\log 2}{\lambda}$ & $\displaystyle \frac{1}{\lambda}$ \\
\\
Beta ${\mathcal B}(\alpha,\beta)$,  $\alpha,\beta \geq 1$ & $\displaystyle \log\left[x^{1-\alpha}(1-x)^{1-\beta}  \right] $ & No expression & --- \\
\\
Gamma $\Gamma(\alpha,\beta)$,  scale $\alpha \geq 1$, shape $\beta > 0$ & $\displaystyle \log\left( x^{1-\alpha} \Gamma (\alpha) \right) +\frac{x}{\beta} +\alpha\log\beta$& No expression & --- \\
\\
Gumbel ${\mathcal G}(\mu,\beta)$, scale $\beta>0$ & $\displaystyle \frac{x-\mu}{\beta} + \log\beta +\exp\left(-\frac{x-\mu}{\beta}\right)$ & $\displaystyle \mu -\beta \log(\log 2)$ & $\displaystyle \frac{\beta}{\log 2}$\\
\\
Weibull ${\mathcal W}(k, \lambda)$, shape $k \geq 1$,\hspace{0.3cm} scale $\lambda > 0$ &$\displaystyle  \log\left(\frac{\lambda}{k}\right) + (1-k)\log\left(\frac{x}{\lambda}\right)+\left(\frac{x}{\lambda}\right)^k $ & $\displaystyle \lambda (\log 2)^{1/k}$ & $\displaystyle \frac{\lambda(\log 2)^{(1-k)/k}}{k}$\\
\hline
  \end{tabular}
    \caption{Standard log-concave probability distributions: $v(\cdot)$ function, median $m$ and Cheeger constant $C_1$ (see Theorem \ref{prop:link2}).}
    \label{tab:logconc}
\end{table}

For general log-concave measures, no analytical expressions are available for the Cheeger constant.
In this latter case or in case of non log-concave but Boltzmann measure, we can estimate the Cheeger constant by numerically evaluating the expression $\sup_{\x \in \mathbb{R}} \frac{\min(F_j(\x), 1-F_j(\x))}{\rho_j(\x)}$.

\section{Numerical tests}\label{sec:test}

\subsection{Derivative sensitivity indices estimates}

A classical estimator for the DGSM is the empirical one and is given below:
   \begin{eqnarray} \label{eq:dsiest}
\widehat{\nu}_{j} &=&  \frac{1}{n} \sum_{i=1}^n \left(\frac{\partial f(\bo{\X}^{(i)})}{\partial \x_j} \right)^2.
\end{eqnarray}
Experimental convergence properties of this estimator are given in Sobol and Kucherenko \cite{sobol09}.

From definition (\ref{eq:decdgsa}), we know that $\displaystyle \frac{\partial f(\bo{\X}^{(i)})}{\partial \x_j}=\frac{\partial g(\bo{\X}^{(i)})}{\partial \x_j}$. Estimator of $D_j^{tot}$ (see equation (\ref{eq:dt})) and estimator (\ref{eq:dsiest}) are based on the same function $g(\cdot)$ and it seems that estimations of these two indices will require approximately the same number of model evaluations in order to converge towards their respective values.

Computation of DGSM and Sobol' indices can be performed with Monte Carlo-like algorithm, such as Latin Hypercube Sampling, quasi-Monte Carlo and Monte Carlo Markov Chain sampling.
Kucherenko $\etal$ \cite{kucherenko09} have shown that quasi-Monte Carlo outperforms Monte Carlo when model has a low effective dimension. Computation of DGSM needs model gradient estimation. For complex models, model gradient computation can easily be obtained by finite difference method. Patelli and Pradlwarter \cite{patelli10} proposed a Monte Carlo estimation of gradient in high dimension. They used an unbiased estimator for gradients and have shown that the number of Monte Carlo evaluations $n \leq d$ is sufficient for gradient computations. In the worst case, their procedure requires the same number of model evaluations than the finite difference method. The method is very efficient when the model has a low effective dimension.  

In the following Sections, we compare the estimates of the Sobol indices ($S_j$ and $S_{T_j}$) and the upper bound of $S_{T_j}$ (see inequality (\ref{eq:ineq})). let denote $\Upsilon_j$, the total sensitivity upper bound:
\begin{equation}\label{eq:upperbound}
\Upsilon_j = C \frac{\nu_j}{D} \;,
\end{equation}
where $D$ is the variance of the model output $f(\bo{\X})$ and $C=4 C_1^2$.
The goal of our numerical tests is just to compare the differences in terms of ranking and not to study the speed of convergence of the estimates.

\subsection{Test on the Morris function}

As a first test, we consider the Morris function (Morris \cite{morris91}) that includes $20$ independent and uniform input variables.
The Morris function is defined by the following equation:
\begin{equation}\label{eq:morris}
 y = \beta_0 + \sum_{i=1}^{20}\beta_i w_i + \sum_{i<j}^{20}\beta_{i,j} w_iw_j + \sum_{i<j<l}^{20}\beta_{i,j,l} w_iw_jw_l + \sum_{i<j<l<s}^{20}\beta_{i,j,l,s} w_iw_jw_lw_s \;,
 \end{equation}
where $\displaystyle w_i=2\left(x_i-\frac{1}{2}\right)$ except for $i=3,5,7$ where $\displaystyle w_i = 2\left(1.1 \frac{x_i}{x_i+1}-\frac{1}{2}\right)$. The coefficient values are:\\
 $\beta_i =20$ for $i=1, 2, \ldots, 10$,\\
 $\beta_{i,j} =-15$ for $i,j=1, 2, \ldots, 6$, $i<j$\\
 $\beta_{i,j,l} =-10$ for $i,j,l=1, 2, \ldots, 5$, $i<j<l$\\
 and $\beta_{1,2,3,4} = 5$.\\
 The remaining first and second order coefficients were generated independently from the normal distribution ${\mathcal N}(0,\,1)$ and the remaining third and fourth coefficient were set to $0$.

We replace the uniform distributions associated with several input variables by different log-concave measures of the Table \ref{tab:logconc} in order to show how the bounds can be used in practical sensitivity analysis.
Table \ref{tab:morris} shows the probability distributions associated to each input of the Morris function.

\begin{table}[!ht]
  \begin{center}
   \begin{tabular}{lclc}
Input & Probability distribution & Input & Probability distribution \\
   \hline
$X1$ & ${\mathcal U}[0,1]$  &      $X11$    &    ${\mathcal U}[0,1]$      \\
$X2$  & ${\mathcal N}(0.5,0.1)$ &  $X12$ &  ${\mathcal N}(0.5,0.1)$     \\
$X3$ & ${\mathcal E}(4)$ &         $X13$  &   ${\mathcal E}(4)$    \\
$X4$ & ${\mathcal G}(0.2,0.2)$ &   $X14$ & ${\mathcal G}(0.2,0.2)$   \\
$X5$  & ${\mathcal W}(2,0.5)$ &    $X15$ &  ${\mathcal W}(2,0.5)$   \\
$X6$ &  ${\mathcal U}[0,1]$ &      $X16$ &      ${\mathcal U}[0,1]$  \\
$X7$ & ${\mathcal U}[0,1]$ &       $X17$ &      ${\mathcal U}[0,1]$  \\
$X8$ & ${\mathcal U}[0,1]$ &       $X18$  &    ${\mathcal U}[0,1]$   \\
$X9$ & ${\mathcal U}[0,1]$  &      $X19$   &   ${\mathcal U}[0,1]$     \\
$X10$ & ${\mathcal U}[0,1]$  &     $X20$  &   ${\mathcal U}[0,1]$    \\
\hline
    \end{tabular}
    \caption{Probability distributions of the input variables of the Morris function}
    \label{tab:morris}
  \end{center}
\end{table}

We have performed some simulations that allow computing the DGSM indices and the Sobol' indices for the $20$ independent factors.
Sobol' indices $S_j$ and $S_{T_j}$ are obtained with the principles described in Saltelli \cite{saltelli02b}, i.e. using two initial Monte Carlo samples of size $10^4$.
For more efficient convergence properties (specially for the case of small indices), the improved formulas proposed by Sobol $\etal$ \cite{sobtar07} for $S_i$ and by Saltelli $\etal$ \cite{saltelli2010b} for $S_{T_i}$ are used.
The approximation errors of these Monte Carlo estimates are calculated by repeating $20$ times the indices estimation and the mean is taken as the estimate. 
With $d=20$ input variables, it leads to $20 \times 10^4 \times (d+2) = 4.4\times 10^6$  model evaluations.
In fact, the size of the Monte Carlo samples have been fitted to achieve acceptable absolute errors (smaller than $1\%$).
However, the objective here is not to compare the algorithmic performances of DGSM and Sobol' indices in terms of computational cost, but just to look at the inputs ranking.

The total Sobol' indices are used in this paper as a reference.
It shows that only the first $10$ inputs have some influence.
Model derivatives are evaluated for each input on a Monte Carlo sample of size $1\times 10^4$ by the finite-difference method (perturbation of $0.01\%$).
Then, DGSM $\nu_j$ require $2.1\times 10^5$ model evaluations.
$\Upsilon_j$ is then computed using equation (\ref{eq:upperbound}) where the variance of the Morris function is estimated to $D=991.521$.
The results are gathered in Table \ref{tab:resmorris}.

\begin{table}[!ht]
  \begin{center}
   \begin{tabular}{l c c c c c c c}
Input &   $S_j$  & $sd$ & $S_{T_j}$ & $sd$ & $\nu_j$ & $C$  & $\Upsilon_j$  \\
\hline
X1 & 0.043 &	0.009 &	0.173 &	0.008 &	2043.820 &	0.101 &	0.209 \\
X2 & 0.007 &	0.003 &	0.029 &	0.002 &	2856.580 &	0.01 &	0.029 \\
X3 & 0.066 &	0.009 &	0.165 &	0.006 &	31653.270 &	0.250 &	7.981 \\
X4 & 0.002 &	0.006 &	0.134 &	0.007 &	2025.950 &	0.333 &	0.680 \\
X5 & 0.035 &	0.005 &	0.055 &	0.003 &	4203.060 &	0.360 &	1.526 \\
X6 & 0.039 &	0.007 &	0.114 &	0.006 &	1337.100 &	0.101 &	0.137 \\
X7 & 0.068 &	0.003 &	0.069 &	0.003 &	6605.960 &	0.101 &	0.675 \\
X8 & 0.156 &	0.007 &	0.157 &	0.007 &	1826.390 &	0.101 &	0.187 \\
X9 & 0.189 &	0.008 &	0.192 &	0.009 &	2249.770 &	0.101 &	0.230 \\
X10 & 0.145 &	0.005 &	0.146 &	0.005 &	1730.400 &	0.101 &	0.177 \\
X11 & 0.000 &	0.001 &	0.002 &	0.001 &	22.630 &	0.101 &	0.002 \\
X12 & 0.000 &	0.000 &	0.000 &	0.000 &	23.940 &	0.01 &	0.000 \\
X13 & 0.000 &	0.001 &	0.001 &	0.000 &	17.670 &	0.250 &	0.004 \\
X14 & 0.001 &	0.001 &	0.003 &	0.001 &	42.850 &	0.333 &	0.014 \\
X15 & 0.000 &	0.001 &	0.001 &	0.001 &	19.870 &	0.360 &	0.007 \\
X16 & 0.000 &	0.001 &	0.002 &	0.001 &	18.860 &	0.101 &	0.002 \\
X17 & 0.000 &	0.001 &	0.002 &	0.001 &	21.400 &	0.101 &	0.002 \\
X18 & 0.000 &	0.001 &	0.002 &	0.001 &	19.950 &	0.101 &	0.002 \\
X19 & 0.000 &	0.001 &	0.004 &	0.001 &	54.380 &	0.101 &	0.006 \\
X20 & 0.000 &	0.001 &	0.004 &	0.001 &	42.250 &	0.101 &	0.004  \\

\hline
    \end{tabular}
    \caption{Sensitivity indices (Sobol' and DGSM) for the Morris function. For the Sobol' indices $S_j$ and $S_{T_j}$, $20$ replicates has been used to get the standard deviation ($sd$).}
    \label{tab:resmorris}
  \end{center}
 \end{table}

In Table \ref{tab:resmorris}, we can first observe that the total sensitivity upper bounds $\Upsilon_j$ are always greater than the total sensitivity indices as expected.
For each input, we distinguish several situations that can occur:
\begin{enumerate}
\item First order and total Sobol' indices are negligible (inputs $X11$ to $X20$).
In this case, we observe that the bound $\Upsilon_j$ is always negligible.
For all the inputs, this test shows the high efficiency of the bound: a negligible bound warrants that the input has no influence.
\item First order and total Sobol' indices significantly differ from zero and have approximately the same value (inputs $X7$ to $X10$).
This means that the input has some influence but no interactions with other inputs.
In this case, the bound $\Upsilon_j$ is relevant (close to $S_{T_j}$), except for $X7$.
The interpretation of the bound gives a useful information about the total influence of the input.
\item First order Sobol' index is negligible while total Sobol' index significantly differs from zero (inputs $X1$ to $X6$).
In this case, the bound $\Upsilon_j$ largely overstimates the total Sobol' index $S_{T_j}$ for $X3$, $X4$ and $X5$.
However, for $X_4$, we have $\Upsilon_4 < 1$ and this coarse information is still usefull.
For the three other inputs, the bound is relevant.
\end{enumerate}

For two inputs ($X3$ and $X5$), results can be judged as strongly unsatisfactory as the bound is useless (larger than $1$ which is the maximal value for a sensitivity index).
We suspect that these results come from:
\begin{itemize}
\item the model non linearity with respect to these inputs (see equation (\ref{eq:morris})),
\item the input distributions (exponential and Weibull).
\end{itemize}
The second explanation seems to be the more convincing as these types of distribution can provide larger values during Monte Carlo simulations.
In this case, departures from the central part of the input domain leads to uncontrolled derivative values of the Morris function.
Indeed, it can be seen that $\nu_j$ is particularly large for $X3$ and $X5$, because of high derivative values in the estimation samples.
Moreover, we have no observed the same results for $X_1$, $X_2$ and $X_4$.

As a conclusion of this first test, we argue that the bound $\Upsilon_j$ is well-suited for a screening purpose.
Moreover, coupling $\Upsilon_j$ interpretation with first order Sobol' indices $S_j$ (estimated at low cost using a smoothing technique or a metamodel, see \cite{salcha00, ioo11}) can bring useful information about the presence or absence of interaction.
For inputs following uniform, normal and exponential distributions,the bound is extremely efficient.
In these particular cases, the  bound is the best one and cannot be improved.


\subsection{A case study: a flood model} \label{sec:case}

To illustrate how the Cheeger constant can be used for factors prioritization, when we use the DGSM, we consider a simple application model that simulates the height of a river compared to the height of a dyke. When the height of a river is over the height of the dyke, flooding occurs. This academic model is used as a pedagogical example in Iooss \cite{ioo11}.
The model is based on a crude simplification of the 1D hydro-dynamical equations of SaintVenant under the assumptions of uniform and constant flowrate and large rectangular sections. It consists of an equation that involves the characteristics of the river stretch:
\begin{equation}
S = Z_v + H -H_d - C_b \quad \mbox{with} \quad H = \left(\frac{Q}{BK_s \sqrt{\frac{Z_m-Z_v}{L} }} \right)^{0.6},
\end{equation}
with $S$ the maximal annual overflow (in meters) and $H$ the maximal annual height of the river (in meters).

The model has $8$ input variables, each one follows a specific probability distribution (see Table \ref{tab:factors}).
Among the input variables of the model, $H_d$ is a design parameter.
The randomness of the other variables is due to their spatio-temporal variability, our ignorance of their true value or some inaccuracies of their estimation. We suppose that the input variables are independent.
\begin{table}[!ht]
  \begin{center}
   \begin{tabular}{lccc}
Input & Description & Unit & Probability distribution \\
   \hline
 $Q$ & Maximal annual flowrate & m$^3$/s & Truncated Gumbel ${\mathcal G}(1013, 558)$ on $[500 , 3000 ]$ \\
 $K_s$ & Strickler coefficient & - & Truncated normal ${\mathcal N}(30, 8)$ on $[15 , +\infty [$ \\
 $Z_v$ & River downstream level & m & Triangular  ${\mathcal T}(49, 50, 51)$ \\
 $Z_m$ & River upstream level  & m  & Triangular  ${\mathcal T}(54, 55, 56)$  \\
 $H_d$ & Dyke height & m &  Uniform ${\mathcal U}[7, 9]$ \\
 $C_b$ & Bank level  & m & Triangular  ${\mathcal T}(55, 55.5, 56)$ \\
 $L$ & Length of the river stretch  & m &  Triangular  ${\mathcal T}(4990, 5000, 5010)$ \\
 $B$ & River width  & m &  Triangular  ${\mathcal T}(295, 300, 305)$ \\
\hline
    \end{tabular}
    \caption{Input variables of the flood model and their probability distributions}
    \label{tab:factors}
  \end{center}
\end{table}

 We also consider another model output: the associated cost (in million euros) of the dyke presence,
\begin{equation}
C_p = \indic_{S>0} +  \left[0.2 + 0.8\left( 1-\exp^{-\frac{1000}{S^4}}\right) \right]\indic_{S \leq 0} + \frac{1}{20}\left(H_d \indic_{H_d>8} + 8 \indic_{H_d \leq 8} \right),
\end{equation}
with $\indic_{A}(x)$ the indicator function which is equal to 1 for $x \in A$ and 0 otherwise.
In this equation, the first term represents the cost due to a flooding ($S>0$) which is 1 million euros, the second term corresponds to the cost of the dyke maintenance ($S \leq 0$) and the third term is the investment cost related to the construction of the dyke. The latter cost is constant for a height of dyke less than $8$ m and is growing proportionally with respect to the dyke height otherwise.

Sobol' indices are estimated with the same algorithms than for the Morris function, using two initial Monte Carlo samples of size $10^5$ and $20$ replicates of the estimates.
It leads to $2\times 10^7$ model evaluations in order to compute first order indices $S_j$ and total indices $S_{T_j}$ (by taking the mean of the $20$ replicates). 
For estimating the DGSM ($\nu_j$, weighted DGSM $\tau_j$ and the total sensitivity upper bound $\Upsilon_j$), a Sobol sequence is used with $1\times 10^4$ model evaluations.

Results of global sensitivity analysis and derivative-based global sensitivity analysis for respectively the overflow $S$ and the cost $C_p$ outputs are listed in Tables \ref{tab:OSI} and \ref{tab:CSI}. Global sensitivity indices show small interaction among input variables for the overflow and the cost outputs. Four input variables ($Q$, $H_d$, $K_s$, $Z_v$) drive the overflow and the cost outputs. This variable classification will serve as reference for comparison issue.

\begin{table}[!ht]
  \begin{center}
   \begin{tabular}{l c c c c c}
Input &   $S_j$  & $S_{T_j} $ & $\nu_j$ & $\tau_j$  & $\Upsilon_j$  \\
\hline
$Q$  & 0.343  & 0.353   & 1.296e-06 &     1.072&  2.807 \\
$K_s$ & 0.130  & 0.139   & 3.286e-03 &     1.033&  0.198 \\
$Z_v$ & 0.185  & 0.186   & 1.123e+00 &  1377.41 &  0.561  \\
$Z_m$ & 0.003  & 0.003  & 2.279e-02 &    33.742&  0.011 \\
$H_d$ & 0.276  & 0.276   & 8.389e-01 &    23.77 &  0.340 \\
$C_b$ & 0.036  & 0.036   & 8.389e-01 &  1268.90 &  0.105 \\
$L$  & 0.000  & 0.000  & 2.147e-08 &     0.268&  0.000 \\
$B$  &  0.000 & 0.000  & 2.386e-05 &     1.070&  0.000 \\
\hline
    \end{tabular}
    \caption{Sensitivity indices for the overflow output of the flood model.}
    \label{tab:OSI}
  \end{center}
 \end{table}

 \begin{table}[!ht]
  \begin{center}
   \begin{tabular}{l c c c c c}
Input &   $S_j$  & $S_{T_j} $ & $\nu_j$ & $\tau_j$  & $\Upsilon_j$  \\
\hline
$Q$  &0.346 & 0.460 & 1.3906e-06 &    2.013    & 3.011e+00 \\
$K_s$ &0.172 & 0.269 & 8.5307e-03 &    1.926    & 5.129e-01 \\
$Z_v$ &0.187 & 0.229 & 1.3891e+00 & 1715.89     & 6.932e-01 \\
$Z_m$ &0.006 & 0.012 & 4.6038e-02 &   68.17     & 2.29e-02 \\
$H_d$ &0.118 & 0.179 & 1.5366e+00 &   44.04     & 6.227e-01 \\
$C_b$ &0.026 & 0.039 & 9.4628e-01 & 1428.69     & 1.180e-01 \\
$L$ & 0.000 & 0.000 & 4.0276e-08 &    0.503    & 2.009e-06 \\
$B$  & 0.001 & 0.001 & 4.4788e-05 &    2.007    & 5.587e-04 \\
\hline
\end{tabular}
    \caption{Sensitivity indices for the cost ouput of the flood model.}
    \label{tab:CSI}
  \end{center}
 \end{table}

Based on derivative sensitivity indices ($ \nu_j$) or weighted derivative sensitivity indices ($\tau_j$) we have obtained another subset of the most influential variables that are $Z_v$, $C_b$, $H_d$, $Z_m$. These results mean that, for example, the maximum annual flowrate ($Q$) does not have any impact on the overflow and the cost output. If we compare these results to the global sensitivity indices, we can infer that they are obviously wrong.
This is easily explained by the fact that the input variables have different unities and that the indices $\nu_j$ and $\tau_j$ have not been renormalized by the constant depending on the probability distribution of $X_j$.

By looking at the total sensitivity upper bound $\Upsilon_j$, the most influential variables are the following: $Q$, $Z_v$, $H_d$, $K_s$ for the overflow output and for the cost output.
It gives the same subset of the most influential variables with some slight differences for the prioritization of the most influential variables.
In conclusion, 
we state that $\Upsilon_j$ can provide correct information on input variance-based sensitivities.

\section{Conclusion} \label{sec:con}

Global sensitivity analysis, that allows exploring numerically complex model and factors fixing setting, requires a large number of model evaluations. Derivative-based global sensitivity method needs a much smaller number of model evaluations (gain factor of $10$ to $100$). The reduction of the number of model evaluations becomes more significant when the model output is controlled by a small number of input variables and when the model does not include much interaction among input variables.
This is often the case in practice.

In this paper, we have produced an inequality linking the total Sobol' index and a derivative-based sensitivity  measure for a large class of probability distributions (Boltzmann measures).
The new sensitivity index $\Upsilon_j$, which is defined as a constant times the crude derivative-based sensitivity, is a maximal bound of the total Sobol' index.
It improves factors fixing setting by using derivative-based sensitivities instead of variance-based sensitivities.

Two numerical tests have confirmed that the bound $\Upsilon_j$ is well-suited for a screening purpose. 
When total Sobol' indices cannot be estimated because of a cpu time expensive model, $\Upsilon_j$ can provide correct information on input sensitivities.
Previous studies have shown that estimating DGSM with a small derivatives' sample (with size from tens to hundreds) allows to detect non influent inputs.
In subsequent works, we propose to use jointly DGSM and first order Sobol' indices.
With these information, an efficient methodology of global sensitivity analysis can be applied and brings useful information about the presence or absence of interaction (see Iooss $\etal$ \cite{ioopop12}).

\section{Acknowlegments}

Part of this work has been backed by French National Research Agency (ANR) through COSINUS
program (project COSTA BRAVA noANR-09-COSI-015). We thank Jean-Claude Fort for helpful discussions and two anonymous reviewers for their valuable comments.

\bibliographystyle{model1b-num-names}
\bibliography{bibi_dgsm}

\end{document}